\documentclass[12pt]{article}
\usepackage{amsfonts}
\usepackage{amsfonts}
\usepackage{graphicx}
\usepackage{amsmath}
\usepackage{amssymb}
\usepackage{amsthm}

\def\LT{{\mathbb{LT}}}

\def\P{\mathop{\hbox{\sf P}}\nolimits}

\def\Z{{\mathbb Z}}

\def\phi{\varphi}

\def\f{\frac}
\def\d{\partial}

\def\B{\Big}
\def\wt{\widetilde}
\def\bs{\backslash}

\newtheorem{theorem}{Theorem}
\newtheorem{lemma}{Lemma}[section]
\newtheorem{definition}{Definition}[section]

\usepackage{graphicx}

\def\figref#1{fig.~\ref{fig.#1}}
\def\eqref#1{(\ref{eq.#1})}

\def\putfigure#1#2{
	\begin{figure}[ht]
	\centering
	\includegraphics{#1}
	\caption{#2}
	\label{fig.#1}
	\end{figure}
}

\def\putfigurescale#1#2{
	\begin{figure}[ht]
	\centering
	\includegraphics[width=4.5in]{#1}
	\caption{#2}
	\label{fig.#1}
	\end{figure}
}

\def\defined#1{{\em #1}}

\def\optional#1{}

\def\proof{\par\noindent{\bf Proof.\ }}
\def\proofof#1{\par\noindent{\bf Proof of #1.\ }}
\def\eop{\vskip 3mm }

\def\calG{{\cal G}}

\def\union{{\cup}}

\def\mucr{{\mu_{cr}}}
\def\wtY{{\widetilde{Y}}}

\DeclareMathAlphabet{\mathpzc}{OT1}{pzc}{m}{it}

\numberwithin{equation}{section}

\theoremstyle{remark}
\newtheorem{rem}{Remark}[section]

\theoremstyle{definition}

\title{Phase transition for the Ising model on the Critical Lorentzian triangulation}
\author{Maxim~Krikun%
  \thanks{Current affiliation: Google Ireland, Dublin, Ireland.}
  \thanks{Partially supported by BQR 2007 UMR7502 IAEM 0039.}$^{~,1}$%
  \and Anatoly~Yambartsev\thanks{Partially supported by the ``Rede
  Matem\'{a}tica Brasil-Fran\c{c}a, CNPq (306092/2007-7)
  and  CNPR ``Edital Universal 2006'' (471925/2006-3).}$^{~,2}$
   }

\begin{document}
\maketitle {\footnotesize

\noindent $^1$
Institut \'Elie Cartan Nancy (IECN), Nancy-Universit\'e, CNRS, INRIA,
Boulevard des Aiguillettes B.P. 239 F-54506 Vand{\oe}uvre-l\`es-Nancy, France. \\
E-mail: krikun@iecn.u-nancy.fr\\
{current affiliation: Google Ireland, Dublin, Ireland}

\noindent $^2$ Department of Statistics, Institute of Mathematics
and Statistics, University of S\~ao Paulo, Rua do Mat\~ao 1010,
CEP 05508--090, S\~ao Paulo SP, Brazil.\\
E-mail: yambar@ime.usp.br}

\begin{abstract}
Ising model without external field on an infinite Lorentzian triangulation
sampled from the uniform distribution is considered. We prove 
uniqueness of the Gibbs measure in the high temperature region
and coexistence of at least two Gibbs measures at low temperature.
The proofs are based on the disagreement percolation method and
on a variant of Peierls method.
The critical temperature is shown to be constant a.s.



\end{abstract}

\section{Introduction}

Triangulations, and planar graphs in general, appear in physics
in the context of 2-dimensional quantum gravity as a model for
the discretized time-space. Perhaps the best understood it the
model of Euclidean Dynamical Triangulations, which can be viewed
as a way of constructing a random graph by gluing together a
large number of equilateral triangles in all possible ways,
with only topological conditions imposed on such gluing.
Putting a spin system on such a random graph can be interpreted
as a {\em coupling of gravity with matter}, and was an object of
persistent interest in physics since the successful application
of matrix integral methods to the Ising model on random lattice
by Kazakov \cite{kazakov}.

More recently, a model of Casual Dynamical Triangulations was
introduced (see \cite{ambjorn} for an overview).
The distinguishing feature of this
model is its lack of isotropy --- the triangulation now has a
distinguished time-like direction, giving it a partial order
structure similar to Minkowski space, and imposing some
non-topological restrictions on the way elementary triangles
are glued.
This last fact destroys the connection between the model and
matrix integrals, in particular the analysis of the Ising model
requires completely different methods (see e.g.~\cite{benedetti}).

From a mathematical perspective, we deal here with nothing but
a spin system on a random graph.
Random graphs, arising from the CDT approach,
were considered in~\cite{myz} under the name of Lorentzian models.
In the present paper we consider the Ising model on such graphs.
When defining the model we pursue the formal Gibbsian approach \cite{georgii};
namely, given a realization of an infinite triangulation,
we consider probability measures on the set of spin configurations
that correspond to a certain formal Hamiltonian.

Our setting is drastically different from e.g.~\cite{kazakov} and~\cite{benedetti}
in that we do not consider ``simultaneous randomness'', when both
the triangulations and spin configurations are included into
one Hamiltonian. Instead we first sample an infinite triangulation
from some natural ``uniform'' measure, and then run an Ising model
on it, thus the resulting semi-direct product measure is ``quenched''.

A modest goal of this paper is to establish a phase transition
for the Ising model in the above described ``quenched'' setting
(the ``annealed'' version of the problem is surely interesting,
but is also more technically challenging, so we don't attempt
it for the moment).
In Section~3 we use a variant of Peierls method to prove
non-uniqueness of the Gibbs measure at low temperature.
Quite surprisingly, proving the uniqueness at high temperature
is not easy -- the difficulty consists in presence of vertices
of arbitrarily large degree, which does not allow for immediate
application of uniqueness criteria such as e.g.~\cite{weitz}.
We resort instead to the method of disagreement percolation~\cite{berg},
and use the idea of ``ungluing'', borrowed from the paper~\cite{bassalygo},
to get rid of vertices of very high degree.

Finally, in Section~5 we show that the critical temperature
is in fact non-random and coincides for a.e. random Lorentzian triangulation.
Section~2 below contains the main definitions
and summarizes some of the results of~\cite{myz}.

We thank E.~Pechersky for numerous useful discussions during the preparation
of this paper.

\section{Definitions and Main Results}

Now we define rooted infinite Lorentzian triangulations in a
cylinder ${\sf C}=S^1\times [0,\infty)$.

\begin{definition}
Consider a connected graph $G$ embedded in a cylinder  ${\sf C}$.
A face is a connected component of ${\sf C} \setminus G$. The face
is a triangle if its boundary meets precisely three edges of the
graph. An embedded triangulation $T$ is such a graph $G$ together
with a subset of the triangular faces of $G$. Let the support
$S(T)\subset {\sf C}$ be the union of $G$ and the triangular faces
in $T$. Two embedded triangulations $T$ and $T'$ are considered
equivalent if there is a homeomorphism of $S(T)$ and $S(T')$ that
corresponds $T$ and $T'$.
\end{definition}
For convenience, we usually abbreviate ``equivalence class of
embedded triangulations'' to ``triangulation''. This should not
cause much confusion. We suppose that the number of the vertices
of $G$ is finite or countable.

\begin{definition}
A triangulation $T$ of ${\sf C}$ is called Lorentzian if the
following conditions hold: each triangular face of $T$ belongs to
some strip $S^1\times[j,j+1], j=0,1,\dots,$ and has all vertices
and exactly one edge on the boundary $(S^1\times\{j\})\cup
(S^1\times\{j+1\})$ of the strip $S^1\times[j,j+1]$; and the
number of edges on $S^1\times\{j\}$ is positive and finite for any
$j=0,1, \dots.$ \end{definition} In this paper we will consider
only the case when the number of edges on the first level $S^1
\times \{0\}$ equal to 1. This is not restriction, only it gives
formulas more clean.

\begin{definition}
A triangulation $T$ is called rooted if it has a root. The root in
the triangulation $T$ consists of a triangle $t$ of $T$, called
the root face, with an ordering on its vertices $(x,y,z)$. The
vertex $x$ is the root vertex and the directed edge $(x,y)$ is the
root edge. The $x$ and $(x,y)$ belong to $S^1\times \{0\}$.
\end{definition}

Note that this definition also means that the homeomorphism in the
definition of the equivalence class respects the root vertex and
the root edge. For convenience, we usually abbreviate
``equivalence class of embedded rooted Lorentzian triangulations''
to ``Lorentzian triangulation'' or LT.

In the same way we also can define a Lorentzian triangulation of a
cylinder ${\sf C}_N=S^1\times [0,N].$ Let $\LT_N$ and
$\LT_{\infty}$ denote the set of Lorentzian triangulations with
support ${\sf C}_N$ and ${\sf C}$ correspondingly.

\subsection{ Gibbs and Uniform Lorentzian triangulations.}
Let $\LT_N$ be the set of all Lorentzian triangulations with only
one (rooted) edge on the root boundary and with $N$ slices. The
number of edges on the upper boundary $S^1\times \{N\}$ is not
fixed. Introduce a Gibbs measure on the (countable) set $\LT_N$:
\begin{equation}
{\sf P} _{N,\mu}(T)= Z_{[0,N]}^{-1} \exp(-\mu F(T)), \label{eq.Gm}%
\end{equation}
where $F(T)$ denotes the number of triangles in a triangulation
$T$ and $Z_{[0,N]}$ is the partition function:%
\[
Z_{[0,N]}=\sum_{T\in\LT_N}\exp(-\mu F(T)).
\]
The measure on the set of infinite triangulations $\LT_\infty$ is
then defined as a weak limit 
\[ \P_\mu := \lim_{N\to\infty} \P_{N,\mu}. \]
It was shown in
\cite{myz} that this limit exists for all $\mu\ge \mucr := \ln2$.
\begin{theorem}[\cite{myz}]\label{Tmyz}
Let $k_n$ be the number of vertices at $n$-th level in a triangulation $T$
for each $n\ge0$.
\begin{itemize}
\item For $\mu>\ln2$ under the limiting measure $\P_\mu$
the sequence $\{ k_n \}$ is a positive recurrent Markov chain.
\item For $\mu = \mucr = \ln 2$ the
sequence $\{k_n\}$ is distributed as the branching process $\xi_n$
with geometric offspring distribution with parameter $1/2$,
conditioned to non-extinction at infinity.
\end{itemize} \end{theorem}
Below we briefly sketch the proof of the second part of Theorem~\ref{Tmyz},
a deeper investigation of related ideas will appear in~\cite{syz}.

Given a triangulation $T\in\LT_N$,
define the subgraph $\tau \subset T$ by taking,
for each vertex $v \in T$, the leftmost edge going from $v$ downwards
(see \figref{lt-tree}).
The graph thus obtained is a spanning forest of $T$,
and moreover, if one associates with each vertex of $\tau$ it's height in $T$
then $T$ can be completely reconstructed knowing $\tau$.
We call $\tau$ the \defined{tree parametrization} of $T$.

\putfigurescale{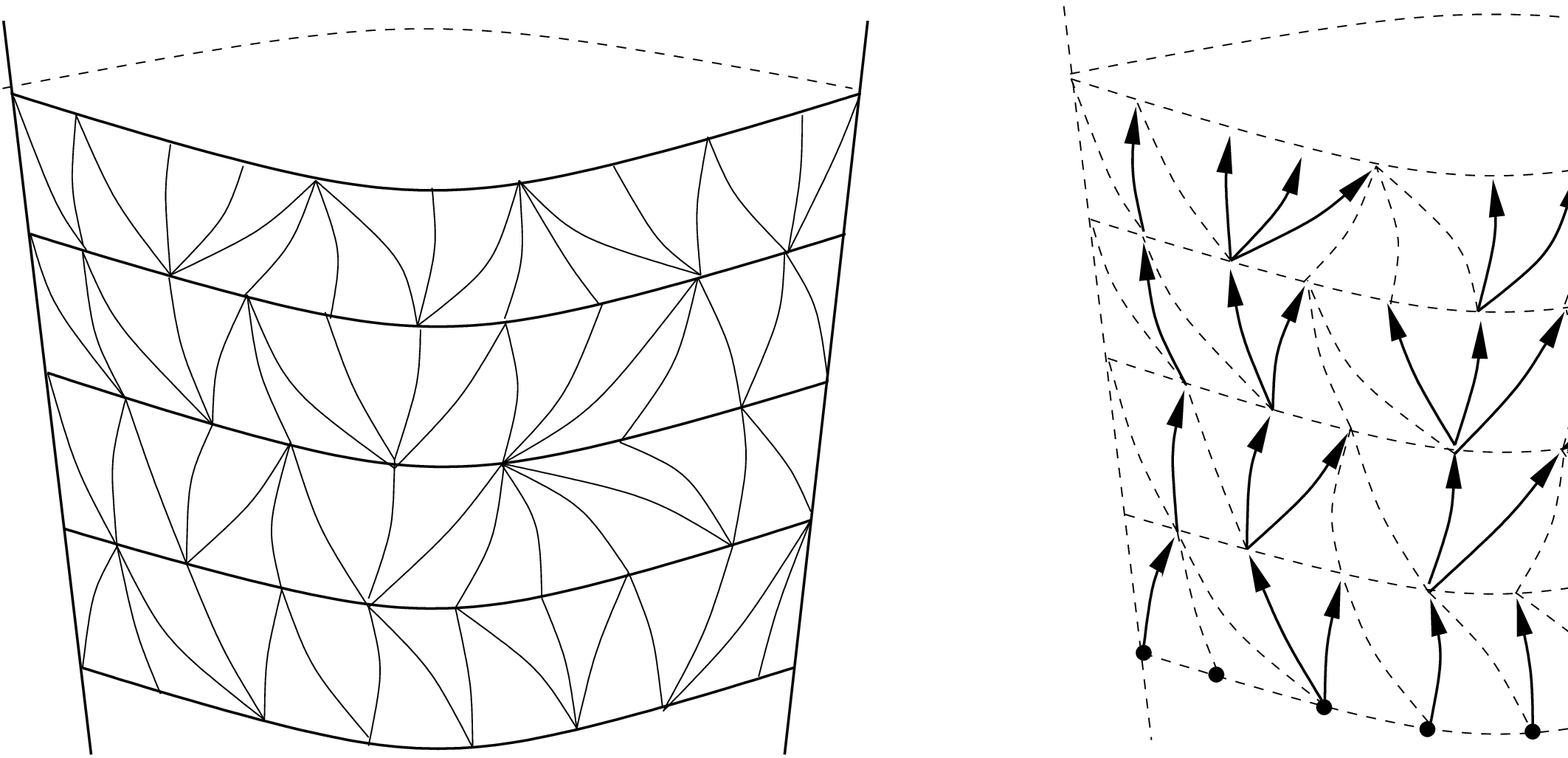}{Tree parametrization}

For every vertex $u\in\tau$ denote by $\delta_u$ it's out-degree,
i.e.~the number of edges of $\tau$ going from $u$ upwards.
Comparing the out-degrees in $\tau$ to the number of vertical
edges in $T$, and comparing the latter to the total number of
triangles, it is not hard to obtain the identity
\begin{equation}\label{eq.sumdeg}
\sum_{u \in \tau \bs S^1\times\{N\}} (\delta_u + 1) = F(T)
\end{equation}
(the sum on the left runs over all vertices of $\tau$ except for the $N$-th level).
Thus under the measure $\P_{N,\mucr}$ the probability of a forest $\tau$
is proportional to
\begin{equation}\label{eq.emucr}
e^{-\mucr F(T) }
 = \prod_{u \in \tau \bs S^1\times\{N\}} \B(\f12\B)^{\delta_u+1}
\end{equation}
which is exactly the probability to observe $\tau$ as a
realization of a branching process with offspring distribution
${\rm Geom}(1/2)$. After normalization we'll obtain, on the left
in \eqref{emucr}, the probability $\P_{N,\mu}(\tau)$ as defined by
\eqref{Gm}, an on the right the {\em conditional} probability to
see $\tau$ as a realization of the branching process $\xi$ given
$\xi_N>0$.
So quite naturally when $N\to\infty$ the
distribution of $\tau$ converges to the Galton-Watson tree,
conditioned to non-extinction at infinity.

In particular it follows from Theorem~\ref{Tmyz} that
\begin{equation} \label{GW}
{\sf P}_\mucr (k_n=m) = Pr(\xi_n = m \mid \xi_\infty >0 ) = m
Pr(\xi_n = m).
\end{equation}

\begin{rem}\label{note1}
The last equality in ($\ref{GW}$) means that the measure $P_\mucr$
on triangulations can be considered as a $Q$-process defined by
Athreya and Ney \cite{AN} for a critical Galton-Watson branching
process. Such a process is exactly a critical Galton-Watson tree
conditioned to survive forever. \end{rem}

We will also note for further use that the offspring generating function
of the branching process $\xi$,
\[ \psi(s) = \sum_{k\ge 0} \B(\f12\B)^{k+1} s^k = \f{s}{2-s}, \]
and the generating function for $\xi_n$ (with initial condition $\xi_0=1$) is
\[ \psi_n(s) = \f{ n-(n-1)s }{ n+1 - ns }. \]

\subsection{Ising model on Uniform Infinite Lorentzian triangulation
-- quenched case.}

Let $T$ be some fixed Lorentzian triangulation, $T\in \LT_\infty.$
Let $T_N$ be the projection of $T$ on the cylinder ${\sf C}_N.$
We associate with every vertex $v$ a spin
$\sigma_v\in \{-1, 1\}$.
Let $\Sigma(T)$ and $\Sigma_N(T)$ denote the set of of all spin
configurations on $T$ and $T_N$, respectively.
The Ising model on $T$ is defined by a formal Hamiltonian
\begin{equation} \label{Hm} H(\sigma) = \sum _{
\langle v, v' \rangle \in V} \sigma_v \sigma_{v'} \end{equation}
where $\langle v, v' \rangle $ means that vertices $v, v'$ are
neighbors, i.e.~are connected by an edge in $T$.
Let $\partial T_N$ be the set of vertices of $T$ that lie on the
circle $S^1 \times \{N+1\}.$ Fix some configuration on the
boundary $\partial T_N$ and denote it $\partial \sigma.$ The Gibbs
distribution with boundary condition $\partial \sigma$ is
defined by the following. Let $V(T_N)$ be the set of all vertices
in $T_N,$ then the energy of configuration $\sigma \in \Sigma
_N(T)$ is \begin{equation} \label{Hb} H_N (\sigma | \partial
\sigma) = \sum_{ \langle v,v' \rangle:\ v,v'\in V(T_N)  } \sigma_v
\sigma_{v'} + \sum _{ \langle v,v' \rangle: v\in V(T_N), v'\in
\partial T_N } \sigma_v \sigma_{v'} \end{equation} which defines
the probability
\begin{equation} \label{GN} P^T_{N, \partial \sigma} ( \sigma ) =
\frac{\exp\{ - \beta H_N (\sigma | \partial \sigma) \} }{ Z_{N,
\partial \sigma} (T)}
\end{equation} where $$ Z_{N, \partial \sigma} (T) = \sum _{\sigma \in
\Sigma_N(T)} \exp\{ - \beta H_N (\sigma | \partial \sigma) \}. $$
When $N\to\infty$, for any sequence of boundary conditions $\d \sigma$,
a limit (at least along some subsequence) of measures $P^T_{N, \d\sigma}$
exists by compactness.
Such a limit is a probability measure on $\Sigma(T)$ with a natural
$\sigma$-algebra, which we refer to as a Gibbs measure.

In general, it is well known that at least one Gibbs measure exists for
the Ising model on any locally finite graph and for any value of
the parameter $\beta$ (see, e.g., \cite{GM} page 71).
It is also known that the existence of more than one Gibbs measure
is increasing in $\beta$,
i.e.~there exists a critical value $\beta_c \in [0,\infty]$ such
that there is a unique Gibbs measure when $\beta > \beta_c$,
and multiple Gibbs measures when $\beta < \beta_c$ (see
\cite{haggstrom} for an overview of relations between percolation
and Ising model on general graphs).

Thus when considering the Ising model on Lorentzian triangulations
it is natural to ask whether the critical temperature is finite
(different from both $0$ and $\infty$), and whether it depends on the
triangulation. In the following two sections we show that the critical
temperature is a.s.~bounded both from $0$ and $\infty$.
In the last section we prove that the critical temperature obeys
a zero-one law and is therefore a.s.~constant.

\section{Phase transition at low temperature}
In this section we prove the following theorem.
\begin{theorem} \label{Th1} There
exists a $\beta_0$ such that for all $\beta \in (\beta_0, \infty)$
there exist at least two Gibbs measures for ${\sf P}_\mucr$-a.e.
$T$.
\end{theorem}

We remind first the classical Peierls method for the Ising model
on $\Z^2$. Let $\Lambda \subset \Z^2$ be a large square box centered at
the origin, and let $P^-_\Lambda$ be the distribution for the
spins in $\Lambda$ under the condition that all of the spins
outside of $\Lambda$ are negative. When the size of the box tends
to infinity, $P^-_\Lambda$ converges to some probability measure
$P^-$, which we call Gibbs measure with negative boundary
conditions. The measure $P^+$ is defined similarly
taking positive boundary conditions.

Due to an obvious symmetry between $P^-$ and $P^+$ we have
\[ P^-(\sigma_0=+1) = P^+(\sigma_0=-1), \]
therefore if the Gibbs measure is unique then necessarily $P^-=P^+$ and
\begin{equation}\label{eq.pleq}
P^-(\sigma_0=+1) = P^-(\sigma_0=-1) = 1/2.
\end{equation}
Thus our goal will be to disprove \eqref{pleq} for large enough $\beta$.

Consider some configuration $\sigma$ such that $\sigma_0=+1$ and
$\sigma_i=-1$ for all $i\notin \Lambda$.
Define the \defined{cluster} $K_+\subset \Z^2$
as the maximal connected subgraph of $\Z^2$,
containing the origin, such that $\sigma_i = +1$ for all $i \in K_+$.
Let $\gamma$ be the contour in the dual graph $(\Z^2)^*$,
corresponding to the outer boundary of $K_+$,
and let $\sigma'$ be the configuration obtained from $\sigma$
by inverting all the spins inside $\gamma$.
Note that $\sigma'_0 = -1$.
Let $C$ be the set of all contours in $(\Z^2)^*$ surrounding the origin,
and consider the application
\begin{equation}
\begin{array}{rrcl}
{\rm \bf inv}: & \Sigma |_{\sigma_0=+1} & \rightarrow & \Sigma |_{\sigma_0=-1} \times C \\
           & \sigma & \rightarrow & (\sigma', \gamma).
\end{array}
\end{equation}
Clearly ${\rm \bf inv}$ is injective: given $(\sigma',\gamma)$ one can
easily reconstruct $\sigma$ by inverting the spins inside
$\gamma$. Also the following property holds:
\begin{equation} \label{eq.plsig}
P^-_\Lambda ( \sigma ) = P^-_\Lambda ( \sigma' ) \cdot e^{-2\beta |\gamma|},
\end{equation}
where $|\gamma|$ denotes the length of $\gamma$. Indeed, for every
edge $\langle i,j \rangle$ traversed by $\gamma$ we have $\sigma_i
\neq \sigma_j$ by construction, but after inversion $\sigma'_i =
\sigma'_j$, so the contribution of the pair $\langle i,j \rangle$
to the Hamiltonian is increased by $2$. On the other hand, for all
other pairs $\langle i,j \rangle$ the contribution to the
Hamiltonian doesn't change. Informally we can write $H(\sigma) -
H(\sigma') = 2|\gamma|$, which is equivalent to~\eqref{plsig}.

Now taking a sum over all $\sigma: \sigma_0=+1$ in \eqref{plsig}
we get
\begin{eqnarray}
P^-_{\Lambda}(\sigma_0=+1)
&=& \sum_{\sigma: \sigma_0 = +1} P^-_\Lambda(\sigma) \nonumber\\
&\le& \sum_{\sigma': \sigma'_0=-1, \gamma \in C}
P^-_\Lambda(\sigma') \cdot e^{-2\beta|\gamma|}
    \nonumber\\
&=& P^-_\Lambda(\sigma_0=-1) \cdot \sum_{n\ge 1} e^{-2\beta n} \# C_n,
\end{eqnarray}
where $C_n$ is the set of contours $\gamma\in C$ of length exactly $n$.

But $\# C_n \le n 3^n$: indeed, every contour $\gamma \in C_n$ intersects the
$x$-axis somewhere between $0$ and $n$.
Starting from this intersection point there are at most $3^n$ distinct
self-avoiding paths of length $n$ in $(\Z^2)^*$,
only a few of them really belonging to $C_n$.
Therefore we have the inequality
\begin{equation}\label{eq.plsum}
P^-_\Lambda(\sigma_0=+1)
\le
P^-_\Lambda(\sigma_0=-1) \cdot \sum_{n\ge 1} e^{-2\beta n} n 3^n,
\end{equation}
and by taking $\beta$ large enough the sum in the right-hand side
of~\eqref{plsum} can be made strictly less than one. Taking the
limit $\Lambda \to \Z^2$ we establish
\[ P^-(\sigma_0=+1) < P^-(\sigma_0=-1). \]
Thus $P^- \neq P^+$ and the Gibbs measure for the Ising model in
$\Z^2$ at inverse temperature $\beta$ is not unique.

\putfigurescale{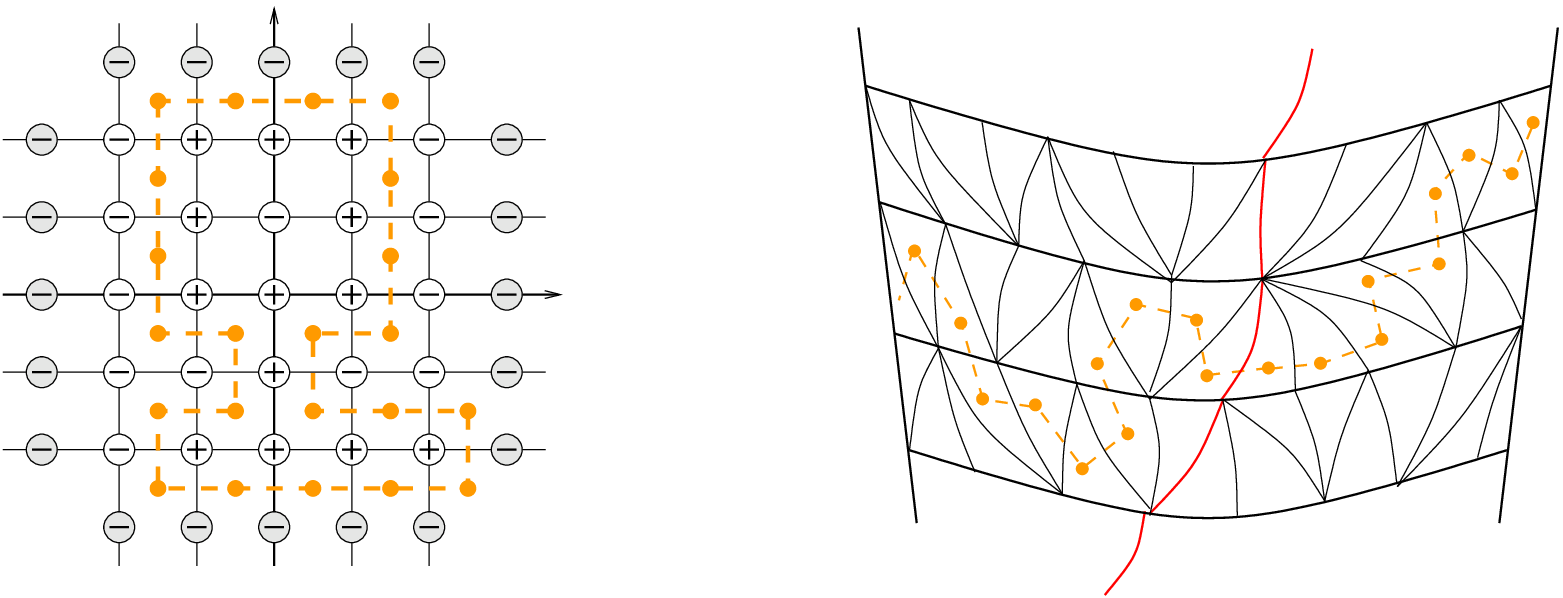}{Peierls' method on $\Z^2$ (left)
and on a random Lorentzian triangulation (right).
Dashed line represents the inversion contour $\gamma$.}

\noindent Next we will slightly generalize the above classical
argument.
\begin{lemma} \label{GenP} Let $G$ be an infinite planar graph
and $G^*$ it's planar dual. Let $v_0$ be a vertex of $G$ and
denote by $C_n(v_0)$ be the set of contours of length $n$ in
$G^*$, separating $v_0$ from the infinite part of the graph.

If the following sum is finite,
\begin{equation}\label{eq.sfin}
\sum_{n\ge 1} \# C_n(v_0) e^{-2\beta n} < \infty,
\end{equation}
then the Gibbs measure for the Ising model on $G$ at inverse
temperature $\beta$ is not unique.
\end{lemma}

\begin{rem} \label{note2}  Strictly speaking, in the above Lemma
    we should also require the graph $G$ to be one-ended;
    i.e.~for any finite subgraph $B\subset G$
    the complement $G\backslash B$ must contain exactly one
    infinite connected component.
    The reason for this is that if the graph $G$ fails to be one-ended
    (e.g.~a doubly-infinite cylinder)
    then the outer boundary of the cluster $K_+$
    may happen to have multiple connected components,
    which slightly complicates the argument.
    We don't consider this case in detail since
    Lorentzian triangulations are one-ended by construction.
\end{rem}

\proof Assume that \eqref{sfin} holds. Then for some large $N$
\[ \sum_{n\ge N} \#C_n(v_0) e^{-2\beta n} < 1. \]
Also, for every $n$ the set $C_n(v_0)$ is finite,
so there exists some large $R$ such that the ball $B_R(v_0)$ contains
all of the $C_n(v_0)$, $n=1\ldots N$.

Let now $P^+$ and $P^-$ be the Gibbs measures for the Ising model on $G$,
constructed
with positive and negative boundary conditions respectively,
and let us compare the events
\[ A_+ := ( \sigma_v = +1 \text{ for all } v \in B_R(v_0) ) \]
and
\[ A_- := ( \sigma_v = -1 \text{ for all } v \in B_R(v_0) ). \]
Proceeding as in the case of $\Z^2$ above, we can show that
\[ P^-(A_+) \le P^-(A_-) \cdot \sum_{n\ge 1} e^{-2\beta n} \# C_n( B_R(v_0)), \]
where $C_n( B_R(v_0))$ now is the set of contours of length $n$ in $G^*$
that surround the whole ball $B_R(v_0)$.
But by construction we have $C_n(B_R(v_0))=\emptyset$ for $n\le N$
and $C_n(B_R(v_0)) \subset C_n(v_0)$ for all other $n$,
therefore
\[ \sum_{n\ge 1} e^{-2\beta n} \# C_n( B_R(v_0))
   \le
   \sum_{n\ge N} e^{-2\beta n} \# C_n(v_0) < 1
 \]
 and
 \[ P^-(A_+) < P^-(A_-). \]
Since by symmetry $P^+(A_-) = P^-(A_+)$, it follows that $P^-\neq
P^+$. \qed

\begin{lemma} \label{l2} Let $T$ be a random Lorentzian triangulation,
and let $v_0$ be the root vertex of $T$. There exists $\beta_0$
such that for every $\beta > \beta_0$ \begin{equation}
\label{finit} \sum_{n\ge 1} \# C_n(v_0) e^{-2 \beta n} < \infty
\quad \mbox{$\P_\mucr$-a.e.},
\end{equation} where $\# C_n(v_0)$ is defined as in Lemma~\ref{GenP}.
\end{lemma}

\proof In order to prove (\ref{finit}) it will be sufficient to
show that the expectation (with respect to the measure $\P_\mucr$)
of the sum is finite
\begin{equation}\label{mcont} {\sf E}_\mucr \sum_{n=0}^\infty \# C_n(0)
e^{-2\beta n} < \infty.
\end{equation}

First of all we choose in any Lorentzian triangulation $T$ a vertical
path $\gamma_\infty = \gamma_\infty(T) = (v_0, v_1, \dots )$
starting at the root $v_0$ and such that $v_i\in S^1\times \{i\}$.
Let $C_{R,n}(v_0) \subset C_n(v_0)$ be the set of contours of
length $n$ which surround $v_0$ and intersect $\gamma_\infty$ at
height $R$. Note that any such contour does not exit from the
strip $[R-n, R+n]$.
Let also $S_{R,n}$ be the number of particles in the tree
parametrization of $T$ at height $R-n$ which have nonempty offspring in
the generation located at height $R+n$. Since every contour from
$C_{R,n}$, in order to surround $v_0$, must cross each of the
$S_{R,n}$ corresponding subtrees, we have
\[ (S_{R,n} > n ) \,\Rightarrow\, (C_{R,n}=\emptyset). \]
On the other hand
\[ C_{R,n} \le 2^n \]
since the contours $C_{R,n}$ live on the dual graph $T^*$, which
has all vertices of degree $3$, thus there are at most $2^n$
self-avoiding paths with a fixed starting point (which is in our
case the intersection with $\gamma_\infty$).

Therefore we have
\begin{equation} \label{est} {\sf E} _\mucr
\sum_{n=0}^\infty \# C_n(v_0) e^{-2\beta n} \le \sum_{n=0} ^\infty
e^{-2\beta n} 2^n \sum_{R=1} ^\infty {\sf P} _\mucr ( S_{R,n} \le
n )
\end{equation}
and for this sum to be finite it's enough to show that sum over
$R$ has polynomial order in $n$.
First let us estimate
$$ \sum_{R=1}
^\infty {\sf P} _\mucr ( S_{R,n} \le n ) \le n^5 + \sum_{R\ge n^5}
^\infty {\sf P} _\mucr ( S_{R,n} \le n ).
$$
and write for some $\epsilon \in (0,1/4)$
\begin{eqnarray} \P _\mucr ( S_{R,n} \le n ) &=& \sum _{k\le R^
\epsilon } {\sf P}_\mucr ( S_{R,n} \le n \mid k_{R-n } = k) {\sf
P} _\mucr ( k_{R-n} =k  ) \nonumber\\
& + & \sum _{k> R^{\epsilon} } {\sf P} _\mucr ( S_{R,n} \le n \mid
k_{R-n } = k ) {\sf P} _\mucr ( k_{R-n} =k  ). \label{eq.sparts}
\end{eqnarray}
The first sum above satisfies the inequality
\begin{equation} \label{est3} \sum _{k\le R^ \epsilon } {\sf P} _\mucr
( k_{R-n} =k ) = {\sf P} _\mucr ( k_{R-n} \le R^ \epsilon ) \le
\Bigl( \frac{ R^\epsilon}{ R-n} \Bigr) ^2.
\end{equation}
Indeed, thanks the representation (\ref{GW}) the generating
function $\phi_n(s)$ of $k_n$ is related to the generating function
$\psi_n(s)$ of the branching process $\xi_n$ described above by
the following equation
\begin{equation} \label{gen.func} \phi_n(s) =
s\psi_n'(s). \end{equation}
Since $\psi_n(s)$ can be calculated explicitly
\begin{equation} \label{gen.GW}
\psi_n(s) = \frac{ n - (n-1) s}{ n+1 - ns},
\end{equation}
from (\ref{gen.func}) and (\ref{gen.GW}) we obtain
\begin{equation}
\phi_n(s) = \frac{ s}{ (1+ n - ns)^ 2},
\end{equation}
\begin{eqnarray} \nonumber
{\sf P} _\mucr ( k_n = k ) &=& \frac{1}{k!} \frac{d^k}{ds^k} \phi_n(s) \Bigl|_{s=0} \\
&=& \f{k n^{k-1}}{(n+1)^k} \le \f{k}{n^2},
\label{prob1}
\end{eqnarray}
and
\[
\sum_{k\le R^\epsilon} \P _\mucr ( k_{R-n}=k ) \le \sum_{k\le
R^\epsilon} \f{k}{(R-n)^2} \le \f{ R^{2\epsilon}}{(R-n)^2},
\]
which proves the estimate (\ref{est3}).
Since $\epsilon < 1/4$ there exists a constant $c_1$ such that
\begin{equation} \sum_{R\ge n^5} \frac{ R^{2\epsilon} } {(R-n)^2}
 < c_1. \end{equation}
The second sum in \eqref{sparts} can be estimated by the
probability $\P _\mucr ( S_{R,n} \le n \mid k_{R-n} = R^\epsilon).$
To do this, let us remind some properties of conditioned or size-biased
Galton-Watson trees (see e.g.~\cite{peres}, \cite{geiger} or \cite{AN}).
If $\tau$ is a Galton-Watson tree of a critical branching process,
conditioned to non-extinction at infinity, then conditionally on
$k_{R-n}$ the subtrees of $\tau$, originating on the $(R-n)$-th level,
are distributed as follows: one subtree, chosen uniformly at random,
has the same distribution as $\tau$ (i.e.~is infinite),
while the remaining subtrees are regular Galton-Watson trees
corresponding to the original branching process (and are finite a.s.).

The probability for a particle of the branching process $\xi$
to survive up to time $2n$ equals $1/(2n+1)$,
thus conditionally on $k_{R-n}$, $S_{R,n}$ is stochastically
minorated by the binomial distribution
with parameters $(k_{R-n}, \f{1}{2n+1})$. 
Using Hoeffding's inequality for binomial distribution we obtain
\begin{equation} \P _\mucr ( S_{R,n} \le n \mid k_{R-n} =
R^\epsilon ) < \exp \Bigl( -2 \frac{ (R^\epsilon/(2n+1) - n)^2
}{R^\epsilon } \Bigr) \end{equation} and the sum over $R>n^5$ is
bounded by some absolute constant $c_2$
\begin{equation} \sum_{R>
n^5} \exp \Bigl( -2 \frac{ (R^\epsilon/(2n+1) - n)^2 }{R^\epsilon
} \Bigr) \le e^2 \sum_{R> n^5} \exp \Bigl( - \frac{ 2R^\epsilon}
{(2n+1)^2 } \Bigr) < c_2. \end{equation}
Thus we have proved that
\begin{equation} \label{est2} {\sf E} _\mucr \sum_{n=0}^\infty C_n(0)
e^{-2\beta n} \le \sum_{n=0} ^\infty e^{-2\beta n} 2^n (n^5 + c_1
+ c_2)
\end{equation} 
and so for all $\beta > \f{\ln 2}{2}$ the sum (\ref{finit}) is
finite. \qed

\noindent
{\it Proof of Theorem~\ref{Th1}.}
Theorem~\ref{Th1} follows immediately from Lemmas~\ref{GenP} and~\ref{l2}.


\section{Uniqueness for high temperature}

In this section we prove the following theorem
\begin{theorem}\label{T2}
There exists a small enough $\beta_h$ such that for every $\beta \in [0, \beta_h)$
for $\P_{\mucr}$-a.e.~Lorentzian triangulation $T$
the Gibbs measure for the Ising model on $T$ is
unique.
\end{theorem}

The proof goes as follows. First, for a fixed triangulation $T$ we
use disagreement percolation to reduce the problem of uniqueness
of a Gibbs measure to the problem of existence of an infinite open
cluster under some specific site-percolation model on $T$. Then we
consider the joint distribution of a triangulation with this
percolation model on it (i.e.~with randomly open/closed vertices),
where the triangulation $T$ is distributed according to
$\P_{\mucr}$, and conditionally on $T$ the state of each vertex in
$T$ is chosen independently (but with vertex-dependent
probabilities). We show that the probability for such a random
triangulation to contain an infinite open cluster is $0$,
therefore this probability is $0$ for
$\P_{\mucr}$-a.e.~triangulation as well, and the theorem follows.

\subsection{Disagreement percolation}

The following result of Van den Berg and Maes \cite{berg}
provides a sufficient condition for the uniqueness of the Gibbs measure,
expressed in terms of some percolation-type problem.
\begin{theorem} 
Let $G$ be a countably infinite locally finite graph with the vertex set $V$,
and let $S$ be a finite spin space.
Let $Y$ be a specification of a Markov field on $G$,
i.e.~for every finite $A\subset V$,
$Y_A(\cdot, \eta)$ is the conditional distribution of the spins on $A$,
given that the spins outside of $A$ coincide with the configuration $\eta \in S^V$.

Define for every $v\in V$
\[ p_v := \max_{\eta, \eta'\in S^V} \, \rho_{\rm var}( Y_{\{v\}}(\cdot, \eta), Y_{\{v\}}(\cdot, \eta' )), \]
where $\rho_{\rm var}$ denotes the distance in variation. Let
$P^{p,G}$ be the product measure under which each vertex $v$ is
open with probability $p_v$ and closed with probability $1-p_v.$

If under site-percolation on $G$, one has
\[
P^{p,G} \{ \mbox{there is an infinite open path}\} = 0,
\]
then the specification $Y$ admits at most one Gibbs measure.
\end{theorem}

In the case of the Ising model, the only spins that affect the
distribution of the spin at a given vertex $v$ are those at
vertices $u$ connected to $v$ by an edge (denoted $u\sim v$). Let
$d_v$ be the degree of $v$ and let $S_v = \sum_{u\sim v}
\sigma_u$. The conditional distribution of $\sigma_v$, given the
values of spins $\sigma_u$ for all $u \sim v$, only depends on
$S_v$ and is given by
\[
P^T\{ \sigma_v = +1 | S_v \} =  \f{ e^{-\beta S_v}}{ e^{-\beta
S_v} + e^{\beta S_v}}, \qquad P^T\{ \sigma_v = -1 | S_v \} =  \f{
e^{\beta S_v}}{ e^{-\beta S_v} + e^{\beta S_v}}.
\]
The distance in variation between two such conditional distributions is maximized
by taking two extreme values of $S_v$, namely $S_v = d_v$ and $S_v = -d_v$,
corresponding to all-plus and all-minus boundary conditions.
Therefore we have
\begin{eqnarray*}
p_v &=&
\f12 \B(
\B| \f{e^{-\beta d_v}}{ e^{-\beta d_v} + e^{\beta d_v}} -
    \f{e^{\beta d_v}}{ e^{-\beta d_v} + e^{\beta d_v}} \B|  \\
    &&{}
+
\B| \f{e^{-\beta d_v}}{ e^{-\beta d_v} + e^{\beta d_v}} -
    \f{e^{\beta d_v}}{ e^{-\beta d_v} + e^{\beta d_v}} \B|
\B)
= \tanh(\beta d_v)
\end{eqnarray*}
and we need to prove the following
\begin{lemma}\label{L4}
Let $T$ be an infinite random Lorentzian triangulation
sampled from the measure $\P_{\mucr}$,
and let each vertex $v\in T$ be open with probability
$p_v = \tanh(\beta d_v)$ and closed otherwise.
Then for small enough $\beta$ the probability
that there is an infinite open path in $T$ is zero.
\end{lemma}

\subsection{Elementary approach to non-percolation}

Imagine for a moment that the graph $G$ that we deal with in
Lemma~\ref{L4} has uniformly bounded degrees. If with positive
probability there exists an infinite open percolation cluster,
then with positive probability this cluster will also be connected
to the origin. Consider the event
\[
A_n = \left\{{
    \mbox{there exists a self-avoiding open path of length $n$}
    \atop\mbox{starting at the origin.}
    }\right\}
\]
The number of paths of length $n$, starting at the origin, is at most $D^n$,
where $D$ is the maximal degree of a vertex in our graph.
Obviously the same estimate is valid for the number of self-avoiding paths.
The probability for a vertex to be open is bounded by $\tanh(\beta D) \le \beta D$,
therefore we can estimate
\[ P^{p,G}(A_n) \le D^n \tanh(\beta D)^n \le \beta^n D^{2n}. \]
Taking $\beta < 1/{D^2}$, we get $P^{p,G} (A_n)\to 0$ as
$n\to\infty$, and therefore the probability that the origin
belongs to an infinite open cluster is zero.

Unfortunately, the above argument can't be applied directly in the context of Lemma~\ref{L4}:
in an infinite random Lorentzian triangulation with probability one one will
encounter regions consisting of vertices of arbitrarily large
degrees, and these regions can be arbitrarily large as well.
On the other hand, we know that in a Lorentzian triangulation the degree
of a typical vertex is not too large (in some sense, which we're not going to make
precise in this paper, it can be approximated
by a sum of two independent geometric distributions plus a constant),
therefore one can hope to show that the vertices of very large degree are so rare
that they don't mess up the picture.

Before proceeding further, let us remark that in the above argument
one can consider, instead of the set of all the self-avoiding paths,
only the paths that are both self-avoiding and \defined{locally geodesic}.
We say that a path $\gamma$ in a graph $G$ is locally geodesic
if whenever two vertices, belonging to $\gamma$, are connected by and edge in $G$,
this edge also belongs to $\gamma$
(in other words, a locally geodesic path avoids making unnecessary detours).
Clearly, if under site percolation on $G$ there exists an infinite
open self-avoiding path starting at $v_0$, there exists also an
infinite self-avoiding locally geodesic path. It is also easy to
see that if $\gamma$ is a locally geodesic path, then every vertex
$v\in\gamma$ has at most two neighbors belonging to $\gamma$. We
will use this observation in the next section.

\subsection{Generalization to random triangulations}

Fix an infinite Lorentzian triangulation $T$ and let $\gamma$ be a
finite self-avoiding, locally geodesic path in $T$, starting at
the root vertex (which we denote by $v_0$). Let $\Gamma\subset T$
be a \defined{$1$-neighborhood} of $\gamma$, i.e.~a
sub-triangulation consisting of the path $\gamma$ together with
all triangles, adjacent to $\gamma$. Such a sub-triangulation
enjoys the \defined{rigidity} property: if $\Gamma$ can be
embedded into
some Lorentzian triangulation, then this embedding is unique%
\footnote{
we require that the embedding maps $v_0$ to the root vertex of $T$,
and sends horizontal/vertical edges of $\Gamma$ to edges of the same type in $T$
}.
We define the event ``$T$ contains $\Gamma$'' (denoted
$\Gamma\subset T$) to be the subset of triangulations in
$\LT_\infty$ for which such an embedding exists.

\begin{lemma}\label{L2B}
Let $\Gamma$ be a $1$-neighborhood of a locally geodesic path $\gamma$,
let $\gamma$ have length $n$ and let the vertices of $\gamma$ have degrees
$d_0, d_1, \ldots, d_n$.

The probability that a random Lorentzian triangulation $T$
contains $\Gamma$ is bounded by
\[
\P_\mucr \{ \Gamma \subset T \}
\le C^n \prod_{i=0}^n \f{1}{d_i^9}
\]
for some absolute constant $C$.
\end{lemma}

We can use this lemma to obtain estimates on percolation probability as follows.
Given the degrees $d_0, \ldots, d_n$ one can completely specify $\Gamma$
using the following collection of numbers
\begin{itemize}
\item for each vertex we need to specify
      how it's degree $d_j$ is split into $d_j^{(up)}$ and $d_j^{(dn)}$,
      i.e.~edges going up and down from the vertex
      --- this gives one integer parameter in the range $[1, d_j]$;
\item also we need to specify the edge along which the path $\gamma$ entered the vertex
      and the edge along which it left; this adds two more parameters.
\end{itemize}
The total number of meaningful combinations for the above
parameters clearly does not exceed $\prod d_j^3$. Now, if $T$ is a
$\P_{\mucr}$-distributed random triangulation with
site-percolation on it, by Lemma~\ref{L2B} the probability for $T$
to contain an open path from the origin with vertex degrees $d_0,
\ldots, d_n$ is bounded by
\[
\prod_j \f{ C }{ d_j^9 } \times \prod_j \tanh(\beta d_j) \times \prod_j d_j^3
\le (\beta C)^n \prod \f{1}{d_j^5}.
\]
Summing over all possible values of $d_0, \ldots, d_n$ we obtain the estimate
\[
\sum_{d_0,\ldots, d_n} (\beta C)^n \prod \f{1}{d_j^5}
\le \B( \beta C \sum_{d\ge 4} \f{1}{d^5} \B)^n
\]
for the probability to have an open path of length $n$ in a random
triangulation. For $\beta$ sufficiently small, the last expression
tends to $0$ as $n\to \infty$, and the probability to have an
infinite open path from the origin is zero.

Note that here we are estimating the annealed probability
(semi-direct product $\P_\mucr\times P^{p,T}$) of the existence of
open path. Clearly, if this annealed probability is zero, then
also for $\P_\mucr$-a.e.~triangulation $T$ the
$P^{p,T}$-probability of an infinite open path is zero as well.

Therefore it remains to prove Lemma~\ref{L2B}.

\subsection{Elementary perturbations}
Consider a finite Lorentzian triangulation $T \in \LT_N$
containing an internal vertex $v$ of very high degree, say~$d_v
\ge 100$. Let $d_v^{(up)}$ and $d_v^{(dn)}$ be the degree
contributed to $d_v$ by edges going upwards and downwards from
$v$, so that $d_v = d_v^{(up)} + d_v^{(dn)} + 2$. Let $O(v)$ be a
$1$-neighborhood of $v$, i.e.~all the triangles adjacent to $v$,
and let $T'$ be the complement to $O(v)$ in $T$.

Let us estimate $\P_{N,\mu}(T | T')$, i.e.~the probability to see
$T$ as a realization of a random Lorentzian triangulation,
conditionally on the event that this triangulation contains $T'$.
For this purpose consider a modification of the $1$-neighborhood
of $v$, which consists in adding a pair of triangles as shown on
\figref{ins1}.
\putfigure{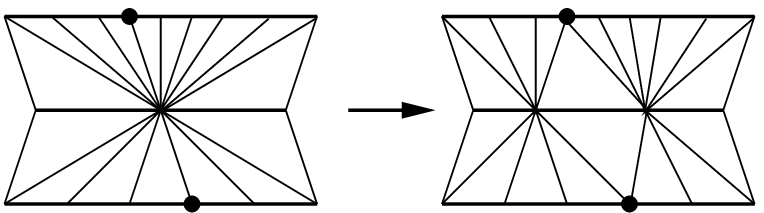}{Inserting a pair of triangles at a high-degree
vertex $v$. The position to insert new triangles is determined by
the choice of two neighboring vertices on the layers below and
above $v$.}
There are $d_v^{(up)} \times d_v^{(dn)}$ to distinct ways to make such a modification,
and every time the Hamiltonian of the Gibbs measure~\eqref{Gm}
is increased by $2\mu$ (since there are $2$ triangles added)

Therefore we can estimate
\[ \P_{N,\mu} (T \,|\, T' ) \le  \f{e^{2\mu}}{ d_v^{(up)} d_v^{(dn)}}
\le \f{ e^{2\mu}}{ (d_v-2)/2}. \]

A better estimate can be obtained by applying $k$
elementary perturbations at $v$ simultaneously.
Such a set of perturbations is specified by a sequence of pairs
$\{(u^{(up)}_i, u^{(dn)}_i)\}_{i=1}^k$ of vertices
lying on the layers above and below $v$ (\figref{ins2}).

\putfigure{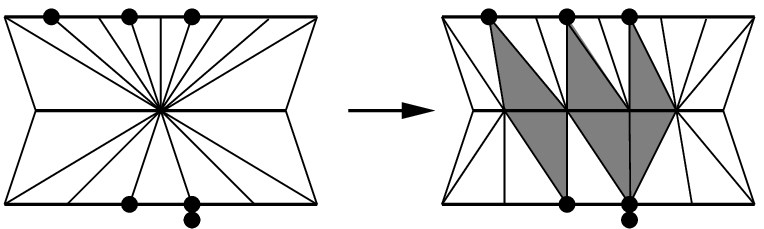}{Multiple insertions at the same vertex.
    Note that one of the edges is used twice for insertions,
    and the corresponding vertex has two dots.}

In order to assure that all the perturbations can be applied unambiguously,
we require that the vertices $u_1^{(up)}, \ldots u_k^{(up)}$ to be
ordered from left to right, as well as $u_1^{(dn)}, \ldots u_k^{(dn)}$.
We may however allow some vertices to coincide -- otherwise it would
be impossible to apply $n$ perturbations at once in a vertex whose
up-degree or down-degree are not large enough.

If $d_v$ is large, one of $d_v^{(up)}$ or $d_v^{(dn)}$ is at least
$(d_v -2)/2$, and therefore the number of distinct ways to choose,
say, $k=10$ modifications is at least
\[ { (d_v -2)/2 \choose 10 }  = C_{10} d_v^{10} + O(d_v ^9). \]

Now consider a path $\gamma \subset \Gamma$ with vertices of degrees $d_0, \ldots, d_n$.
Let $W_\Gamma$ be the set of possible modifications we can make to $\Gamma$,
by inserting exactly $10$ pairs of triangles at vertices of high degree
(say, $d_i \ge 100$), and leaving all other vertices intact.
Then for some constant $\widetilde C_{10}$
\[ \# W_\Gamma \ge \prod \widetilde C_{10} d_j^{10}. \]

Let $\LT_{\Gamma,N} \subset \LT_N$ be the set of Lorentzian triangulations that contain~$\Gamma$.
Clearly, each modification $w\in W_\Gamma$ can be applied to every $T \in \LT_{\Gamma,N}$,
producing a modified triangulation which we denote by $w(T)$.
Consider the application of modifications $w$ to a triangulation $T$
as a mapping
\begin{equation}\label{eq.app}
\begin{array}{rccl}
{\rm \bf app}: & \LT_{\Gamma,N} \times W_\Gamma & \rightarrow & \LT_N \\
           & (T, w) & \rightarrow & w(T).
\end{array}
\end{equation}

Since any modification to $T \in \LT_{\Gamma,N}$
adds at most $20 n$ new triangles,
we have, for all $(T,w) \in \LT_\Gamma \times W_\Gamma$
\[ \P_{N,\mu}( w(T) ) \ge e^{-20n \mu} \P_{N,\mu}(T). \]
Taking a sum over all $(T,w)$ we get
\begin{eqnarray}\label{eq.twsum}
\sum_{(T,w)\in \LT_{\Gamma,N} \times W_\Gamma} \P_{N,\mu}(w(T))
&\ge& e^{-20n \mu} \# W_\Gamma \sum_{T\in \LT_{\Gamma,N}}
\P_{N,\mu}(T) \\ && {} = e^{-20n \mu} \# W_\Gamma
\P_{N,\mu}(\Gamma). \nonumber
\end{eqnarray}

\subsection{Overcounting via random reconstruction}

In order to get an estimate on $\P_{N,\mu}(\Gamma)$ from~\eqref{twsum},
it would be helpful if the application~\eqref{app} was an injection.
But this is not easy to prove (if at all true!),
so instead we will estimate the \defined{overcounting}
\[
{\rm over}_\Gamma(T')
:= \#\{ (T,w) \in \LT_{\Gamma,N} \times W_\Gamma \,|\, w(T) = T' \}.
\]

\begin{lemma}\label{L2C}
Given a path $\gamma$ with $1$-neighborhood $\Gamma$,
with vertex degrees $d_0, \ldots, d_n$, we have for any $T'\in \LT_N$
\begin{equation}\label{eq.overt}
{\rm over}_\Gamma(T') \le \prod 12 (d_j + 10).
\end{equation}
\end{lemma}

\proof
First note that each elementary perturbation (i.e.~insertion of two adjacent
triangles) can be undone by collapsing the newly inserted horizontal edge.
Therefore in order to reconstruct the pair $(T,w)$, knowing $T'=w(T)$,
it will be sufficient to identify the edges of $T'$ that were added by $w$.

It's not clear whether we can perform such a reconstruction deterministically,
but we certainly can achieve it {\em with some probability} with the help
of the following randomized algorithm:
\begin{itemize}
\item{} the algorithm starts from the root vertex
        and at each step moves to an adjacent
        vertex, chosen uniformly at random;
\item{} it also modifies the triangulation along his way as following
\item{} each time arriving to a new vertex
        it should decide whether to contract a sequence of $10$ horizontal edges
        (at least one of these edges must be adjacent to the current location),
        or do nothing.
        This gives $12 = 11 + 1$ possibilities,
        among which one is chosen uniformly at random.
\item{} after $n$ steps, the path traversed by the algorithm is declared to
        be the (conjectured) path $\gamma$, and if it really is,
        and if the resulting triangulation coincides with the original triangulation $t$,
        the reconstruction is considered successful.
\end{itemize}
Most of the time this algorithm will fail,
but with some probability it may actually reconstruct both $T$ and $w$.
Let us now estimate this probability.

First note that the degree of $v_0$ in $T'$ can be larger than the degree of $v_0$
in $T$ because of the modifications made to the second vertex of $\gamma$,
but in any case it will not exceed $d_0 + 10$, since each elementary
modification in $v_1$ increases the degree of $v_0$ at most by $1$.
Therefore, the algorithm will guess the correct direction from $v_0$
with probability at least $1/(d_0 + 10)$.

Assuming the direction it took on the first step was correct,
there is exactly $1$ right choice out of $12$ on the next step
(collapse some $10$ nearby horizontal edges,
 which can be chosen in $11$ different ways, or do nothing),
so the probability not to fail is $1/12$.

Assuming again that
the modifications to the second vertex $v_1$ were undone correctly,
there is at least one edge that will lead
to the next vertex $v_2$.
Since the path $\gamma$ is locally geodesic,
at most two neighbors of $v_1$ belong to $\gamma$,
namely these are $v_0$ and $v_2$.
Modifications at $v_2$ add at most $10$ to the degree of $v_1$,
so with probability at least $1/(d_1 + 10)$
the algorithm will take the correct direction again.

Proceeding this way, we see that the probability to win is at least
\begin{equation}\label{eq.rec}
\prod_{j=0}^n \f{1}{12 (d_j+10)}.
\end{equation}

Now for every $T'$ and for every $(T,w) \in \LT_{\Gamma,N} \times W_\Gamma$
such that $w(T)=T'$ the probability to correctly reconstruct
$(T,w)$ from $T'$ is at least \eqref{rec},
therefore
\[ {\rm over}_\Gamma(T') \le \prod 12 (d_j+10). \]
This finishes the proof of Lemma~\ref{L2C}
\qed

Finally, combining the inequalities \eqref{twsum} and \eqref{overt} we get
\begin{eqnarray*}
e^{-20n \mu} \# W_\Gamma \P_{N,\mu}(\Gamma) &\le& \sum_{(T,w)\in
\LT_{\Gamma,N} \times W_\Gamma} \P_{N,\mu}(w(T))
\\
&\le&
\sum_{T' \in \LT_N} {\rm over}_\Gamma(T') \P_{N,\mu}(T')
\\
&\le&
\prod_{j=0}^n 12 (d_j+10) {\sum_{T' \in \LT} \P_{N,\mu}(T') }
\\
&\le&
 \prod_{j=0}^n 12 (d_j+10),
\end{eqnarray*}
and
\[
\P_{N,\mu}(\Gamma)
\le
\f{ e^{20n \mu} \prod\limits_{j=0}^n 12 (d_j+10) }
  { \prod\limits_{j=0}^n \widetilde C_{10} d_j^{10} }
\le
\prod_{j=0}^n \f{e^{20\mu}{\widetilde C} }{ d_j^9 }.
\]
Since the last estimate is uniform in $N$,
it holds also for the limiting measure~$\P_\mu$,
and setting $\mu=\mucr$ we obtain the statement of Lemma~\ref{L2B}.
This also finishes the proof of Theorem~\ref{T2}.



\section{Critical temperature is constant a.s.}
Consider the critical temperature $\beta_c(G)$ of the Ising model on a graph $G$
as a function of $G$. In the above two sections we have shown that
when $T$ is a $P_{\mucr}$-random Lorentzian triangulation,
we have $\beta_c(T) \in [\beta_h, \beta_0]$ a.s.
In this section we show that in fact
\begin{theorem}\label{T3}
$\beta_c(T)$ is constant $\P_\mucr$-a.s.
\end{theorem}
The proof relies on two lemmas.
\begin{lemma}\label{Lsimplex}
Let $G$, $G'$ be two locally finite infinite graphs
that differ only by a finite subgraph,
i.e.~there exist finite subgraphs $H\subset G$, $H'\subset G'$
such that $G\bs H$ is isomorphic to $G'\bs H'$.
Then $\beta_c(G) = \beta_c(G')$.
\end{lemma}
\proof
This statement is widely known and follows from Theorem~7.33~in~\cite{georgii}.

More exactly,
consider the Ising model on $G$ at some fixed inverse temperature $\beta\in(0,\infty)$.
Let $Y$ be the corresponding specification (collection of conditional distributions)
and denote by $\calG(Y)$ the set of Gibbs measures, satisfying this specification.
It is an easy fact that $\calG(Y)$ is a convex set in the linear space of
finite measures on the set of configurations.

Any specification $Y$ on $G$ can be naturally restricted
to a specification $\wtY$ on~$G\bs H$.
Namely,
for each finite $A \subset G\bs H$ and each $\eta \in \{-1,1\}^{G\bs H}$ put
\begin{eqnarray}
\wtY_A(\cdot, \eta) &=&
\P\{ \sigma_i = (\cdot)_i\,\forall i \in A;\,\,
     \sigma_i = \eta_i\,\forall i\in (G\bs H)\bs A
     \}.
     \nonumber\\
&=& \sum_{\wt\eta \in \{-1,+1\}^H} Y_{A\union H}( (\cdot)\union \wt\eta \,|\, \eta )
\label{eq.wyah}
\end{eqnarray}
(the last sum runs over all spin configurations in $H$).
Informally, the restricted specification $\wtY$ can be interpreted as
an Ising model on $G$, but with the spins in $H$ being hidden from the
observer.

%
Let now $Y'$ be the specification for the Ising model on $G'$,
and let $\wtY'$ be an analogous restriction of $Y'$ to $G'\bs H'$.
Since $G\bs H$ and $G'\bs H'$ are isomorphic, the specifications
$\wtY$ and $\wtY'$ are defined over the same graph,
and it follows from \eqref{wyah} that they are \defined{equivalent}
in that there exists a constant $c>1$ such that
\[ \f1c <  \f{ \wtY_A(\cdot \,|\, \eta) }{ \wtY'_A (\cdot \,|\, \eta) } < c \]
for all possible values of $A$, $\eta$ and $(\cdot)$.
Theorem 7.33 in \cite{GM} then implies that
$\calG(\wtY)$ is affinely isomorphic to $\calG(\wtY')$.
In particular $\calG(\wtY)$ consists of a single element
-- a unique Gibbs measure -- if and only if $\calG(\wtY')$ does.

Since any Gibbs measure on spin configurations in $G\bs H$,
satisfying $\wtY$, can be uniquely extended to a Gibbs measure on $G$, satisfying $Y$,
and the same holds for $G'$, $H'$, $Y'$, the lemma follows.
\eop

\begin{lemma}\label{Lspinal}
A Galton-Watson tree for a critical branching process,
conditioned to non-extinction at infinity,
can be described as following
\begin{itemize}
\item{} it contains a single infinite path (a \defined{spine}), $\{v_0, v_1, \ldots\}$,
        starting at the root vertex;
\item{} at each vertex $v_i$ of the spine a pair of finite trees $(L_i, R_i)$
        is attached, one of each side of the spine;
\item{} and the pairs $(L_i, R_i)$ are independent and identically distributed.
\end{itemize}
\end{lemma}

\begin{rem} In fact, the law of $(L_i, R_i)$ can be explicitly expressed
in terms of the original branching process, but for our purposes
it's enough to know that the trees are i.i.d.
\end{rem}

\proof See e.g.~\cite{peres} or \cite{geiger}.
\eop

\proofof{Theorem~\ref{T3}}
Consider now the tree parametrization $\tau$ of a $\P_\mucr$-distributed
infinite Lorentzian triangulation $T$.
By Lemma~\ref{Lspinal}, $\tau$ is encoded by
a sequence of pairs of finite trees $(L_i,R_i)_{i\in\Z}$.
Clearly, this is an exchangeable sequence.
Let $\sigma_{i j}$ be a permutation of the $i$-th and $j$-th pairs;
then $\sigma_{i j}\tau$ is a tree-parametrization of some other
infinite Lorentzian triangulation $\sigma_{i j}T$,
which only differs from $T$ by a finite subgraph.
By Lemma~\ref{Lsimplex} $\beta_c(T) = \beta_c(\sigma_{i j} T)$,
and Hewitt-Savage zero-one law applies to $\beta_c$, considered
as a function of the sequence $(L_i,R_i)_{i\in\Z}$.
Thus $\beta_c(T)$ is constant a.s.
\eop

\bibliographystyle{amsplain}
\bibliography{lt}

\providecommand{\bysame}{\leavevmode\hbox to3em{\hrulefill}\thinspace}
\providecommand{\MR}{\relax\ifhmode\unskip\space\fi MR }
\providecommand{\MRhref}[2]{%
  \href{http://www.ams.org/mathscinet-getitem?mr=#1}{#2}
}
\providecommand{\href}[2]{#2}
\begin{thebibliography}{10}

\bibitem{AN}
Krishna~B. Athreya and Peter~E. Ney, \emph{Branching processes},
  Springer-Verlag, New York, 1972, Die Grundlehren der mathematischen
  Wissenschaften, Band 196. \MR{MR0373040 (51 \#9242)}

\bibitem{bassalygo}
L.~A. Bassalygo and R.~L. Dobrushin, \emph{Uniqueness of a {G}ibbs field with a
  random potential---an elementary approach}, Teor. Veroyatnost. i Primenen.
  \textbf{31} (1986), no.~4, 651--670. \MR{MR881577 (88i:60160)}

\bibitem{benedetti}
D.~Benedetti and R.~Loll, \emph{Quantum gravity and matter: Counting graphs on
  causal dynamical triangulations}, General Relativity and Gravitation
  \textbf{39} (2007), 863--898.

\bibitem{geiger}
Jochen Geiger, \emph{Elementary new proofs of classical limit theorems for
  {G}alton-{W}atson processes}, J. Appl. Probab. \textbf{36} (1999), no.~2,
  301--309. \MR{MR1724856 (2001k:60119)}

\bibitem{georgii}
Hans-Otto Georgii, \emph{Gibbs measures and phase transitions}, de Gruyter
  Studies in Mathematics, vol.~9, Walter de Gruyter \& Co., Berlin, 1988.
  \MR{MR956646 (89k:82010)}

\bibitem{GM}
Hans-Otto Georgii, Olle H{\"a}ggstr{\"o}m, and Christian Maes, \emph{The random
  geometry of equilibrium phases}, Phase transitions and critical phenomena,
  Vol. 18, Phase Transit. Crit. Phenom., vol.~18, Academic Press, San Diego,
  CA, 2001, pp.~1--142. \MR{MR2014387 (2004h:82022)}

\bibitem{haggstrom}
Olle H{\"a}ggstr{\"o}m, \emph{Markov random fields and percolation on general
  graphs}, Adv. in Appl. Probab. \textbf{32} (2000), no.~1, 39--66.
  \MR{MR1765172 (2001g:60246)}

\bibitem{kazakov}
V.A. Kazakov, \emph{Ising model on a dynamical planar random lattice: Exact
  solution}, Phys. Lett. A \textbf{119} (1986), 140--144.

\bibitem{ambjorn}
R.~Loll, J.~Ambjorn, and J.~Jurkiewicz, \emph{The universe from scratch},
  Contemporary Physics \textbf{47} (2006), 103--117.

\bibitem{peres}
Russell Lyons, Robin Pemantle, and Yuval Peres, \emph{Conceptual proofs of
  {$L\log L$} criteria for mean behavior of branching processes}, Ann. Probab.
  \textbf{23} (1995), no.~3, 1125--1138. \MR{MR1349164 (96m:60194)}

\bibitem{myz}
V.~Malyshev, A.~Yambartsev, and A.~Zamyatin, \emph{Two-dimensional {L}orentzian
  models}, Mosc. Math. J. \textbf{1} (2001), no.~3, 439--456, 472.
  \MR{MR1877603 (2002j:82055)}

\bibitem{syz}
V.~Sisko, A.~Yambartsev, and A.~Zamyatin, \emph{Uniform infinite lorentzian
  triangulation and critical branching process}, in preparation.

\bibitem{berg}
J.~van~den Berg and C.~Maes, \emph{Disagreement percolation in the study of
  {M}arkov fields}, Ann. Probab. \textbf{22} (1994), no.~2, 749--763.
  \MR{MR1288130 (95h:60154)}

\bibitem{weitz}
Dror Weitz, \emph{Combinatorial criteria for uniqueness of {G}ibbs measures},
  Random Structures Algorithms \textbf{27} (2005), no.~4, 445--475.
  \MR{MR2178257 (2006k:82036)}

\end{thebibliography}
\end{document}